\documentclass[11pt]{article}

\usepackage[linewidth=1pt]{mdframed}
\usepackage[bookmarks]{hyperref}
\usepackage{pdfpages}
\usepackage{epsfig}
\usepackage{url}
\usepackage{graphicx,color}
\usepackage[psamsfonts]{amssymb}
\usepackage{amsmath}
\usepackage{indentfirst}
\usepackage{mathrsfs}
\usepackage{booktabs}
\usepackage{tabularx}
\usepackage{setspace}
\usepackage{color}

\usepackage{comment}
\usepackage{bm}
\usepackage{bbm}
\usepackage{dirtytalk}
\usepackage{epstopdf}
\usepackage{mathtools}
\DeclarePairedDelimiter{\floor}{\lfloor}{\rfloor}

\newcommand{\Pm}{{\mathbb P}}
\newcommand{\F}{{\mathcal{F}}}
\newcommand{\B}{{\mathcal{B}}}
\newcommand{\R}{{\mathbb R}}
\newcommand{\N}{{\mathbb N}}
\newcommand{\Z}{{\mathbb Z}}
\newcommand{\E}{{\mathbb E}}
\newcommand{\ls}{{\;\leqslant\;}}
\usepackage{environ}
\NewEnviron{as}{%
  \begin{align}\begin{split}
    \BODY
  \end{split}\end{align}}
  
\usepackage{ntheorem}
\theoremseparator{.}
\newtheorem*{conjecture}{Conjecture}
\newtheorem{definition}{Definition}

\DeclareMathOperator{\err}{err}
\DeclareMathOperator*{\argmax}{arg\,max}

\title{The ARMA Point Process and its Estimation}
\author{Spencer Wheatley$^{1*}$, Michael Schatz$^{1*}$, and Didier Sornette$^1$\\
$^1$ {\small ETH Zurich, Department of Management, Technology and Economics, Switzerland}\\
{\small e-mails: swheatley@ethz.ch, mschatz@ethz.ch, and dsornette@ethz.ch}\\
{\small *These authors contributed equally.}}

\begin{document}
\maketitle

\begin{abstract}
We introduce the ARMA (autoregressive-moving-average) point process, which is a Hawkes process driven by a Neyman-Scott process with Poisson immigration. It contains both the Hawkes and Neyman-Scott process as special cases and naturally combines self-exciting and shot-noise cluster mechanisms, useful in a variety of applications.
The name ARMA is used because the ARMA point process is an appropriate analogue of the ARMA time series model for integer-valued series. As such, the ARMA point process framework accommodates a flexible family of models sharing methodological and mathematical similarities with ARMA time series. We derive an estimation procedure for ARMA point processes, as well as the integer ARMA models, based on an MCEM (Monte Carlo Expectation Maximization) algorithm. This powerful framework for estimation accommodates trends in immigration, multiple parametric specifications of excitement functions, as well as cases where marks and immigrants are not observed.
\end{abstract}

\pagebreak
\clearpage

\section{Introduction}

Mixtures of exo (exogenous) and endo (endogenous) processes pervade science and nature \cite{Sornetteendoexo05}. For instance, in seismology there are main-shocks and after-shocks \cite{helmstetter2002subcritical}. In epidemiology, basic incidence of diseases are followed by contagious outbreaks. In financial markets, there are events that lack (identifiable) precursors followed as well as positive feedbacks and self-excited trading activity \cite{filimonov2012quantifying,HardimanBouchaud2013,Muzy2013hawkes,bacry2015hawkes}. Getting the endo-exo distinction right, as well as the endo mechanism, is fundamental to understanding such processes. In practice, in the absence of physical theory, one typically has a sample of points in time and/or space, which exhibit some kind of clustering. Point-process models featuring a specific generating process can be considered, often with a causal structure and interpretation. Both temporal and spatial versions of these processes exist.

A pre-dominant model type for this purpose is based on the Hawkes process. This is largely due to its parsimonious representation of exo and endo activity, and ease of inference within the likelihood framework. The Hawkes process features exogenous activity through Poissonian immigration, and endogenous activity through an auto-regressive self-excitation. In fact, an equivalence between the Hawkes process and the INAR (integer valued autoregressive) time series model has been proven \cite{kirchner2016hawkes}. A second type of model, the so-called shot noise process, features single bursts of points and has been the scope of insurance risk-theoretic analysis \cite{dassios2003pricing,albrecher2006ruin,kluppelberg1995explosive}. This process can be understood as the analogue of the INMA (integer valued moving average) time series model. More general parallels between time series and point processes have been explored \cite{brillinger1994time}.

It is therefore natural to identify the point process analogue of the INARMA (integer valued ARMA) time series model \cite{Fokianos,McKenzie2003,weiss2008thinning,weiss2008serial}. Here we introduce the ARMA (autoregressive-moving-average) point process and argue that it is in fact such an analogue. This extends the Hawkes process to also be driven by shot-noise type bursts -- the point process analogue to a MA time series -- in addition to the auto-regressive self-excitation. This process is closely related to the \emph{dynamic contagion} point process \cite{dassios2011dynamic}, while a more general class of Hawkes processes with general immigration has been introduced in \cite{Bremaud2002} and further analyzed in \cite{Boumezoued2016}. 
An example for the ARMA point process, is in the analysis of high frequency price changes in financial markets, often studied with Hawkes process \cite{bacry2015hawkes,chavez2012high,filimonov2012quantifying}. The ARMA is an improved model, allowing for exogenous shocks and clustering due to diverse mechanisms, such as order splitting and the near-simultaneous independent actions of multiple market participants in response to an exogenous shock. Other applications, as argued in \cite{dassios2011dynamic}, include credit default claims and insurance claims, both being subject to exogenous shocks as well as endogenous contagion.

Statistical estimation of models with a moving average (shot noise) component based on
MLE (maximum likelihood estimation) is complicated by the fact that the innovations are not observed, so other less efficient or approximate methods are relied upon \cite{brannas2001estimation,tanaka2008parameter,prokevsova2013asymptotic,dassios2005kalman,brillinger2012statistical}. To overcome this difficulty and enable MLE for the ARMA point process, we derive an MCEM (Monte-Carlo expectation-maximization) algorithm, extended from the spatial point process literature on shot noise process  \cite{Moller2003,moller2003statistical} to our setting with both self-exiting and shot noise components. While, as noted above, similar processes have been studied in the literature \cite{dassios2011dynamic,Bremaud2002,Boumezoued2016}, we are not aware of an estimation procedure based on MLE that allows one to
\begin{enumerate}
\item include time-dependence (trends) in immigration,
\item include general parametric or non-parametric excitation functions, and 
\item deal with unobserved immigrants and marks.
\end{enumerate}
In simulation studies, this MCEM algorithm is shown to perform well. However, it is found that estimation of the true log-likelihood is deeply inefficient. This prohibits model selection on a likelihood basis and leaves this desirable feature as an open problem. 
A similar approach using MCEM for the estimation of INARMA time series is discussed in section 4 of \cite{enciso2009efficient} but found to be feasible only for INAR submodels (self-exciting activity). 

Below we first define the ARMA point process, derive basic properties and discuss connections to closely related models from the point process literature (section \ref{sec:ARMA}). Then we establish a formal analogy to the INARMA time series (section \ref{sec:inarma}) and carefully develop the estimation procedure based on Monte-Carlo expectation maximization (section \ref{sec:EM}). We close with a simulation study to demonstrate performance of the estimation procedure and its ability to select correct submodels within the ARMA framework (section \ref{sec:simstudy}).

\pagebreak
\clearpage

\section{The ARMA Point Process}\label{sec:ARMA}
Below we will work with point processes defined as random measures on $\R$. For such a process $N$ and a Borel set $A\subseteq\R$, let $N(A)$ denote the (random) number of points in $A$. In the special case of a half open interval $(a,b]$, we will write $N((a,b])=N(a,b)$. We will use $N_t=N(0,t)$ for $t\in[0,\infty)$ and $N_t=-N([t,0))$ for $t\in(-\infty,0)$. Sometimes it is useful to view the process on a finite window $[0,t]$ as the induced random sequence $T_t:\Omega\rightarrow \R^{\N}$ consisting of all the times $s\in[0,t]$ where $\Delta N_s\neq 0$. A single realisation of $T$ will then be a vector denoted by $T_t(\omega)=\boldsymbol{t}_{1:n}=(t_1,\dots,t_n)$ for $n=N_t(\omega)$ and some (fixed) $\omega\in\Omega$. For a process $N$ defined on a probability space $(\Omega,\F,\Pm)$, we will denote by $\sigma(N)$ the history of the process $N$ given by $\sigma(N)=(\F_t)_{t\in(-\infty,\infty)}$, where $\F_t\subseteq\F$ is the $\sigma$-algebra generated by the evolution of $N$ up to and including $t$, that is, $\F_t=\sigma(N_s,s\leq t)$.

\subsection{Setting and Definition}\label{sec:def}
\begin{figure}[h!]
\centerline{\includegraphics[width=10cm]{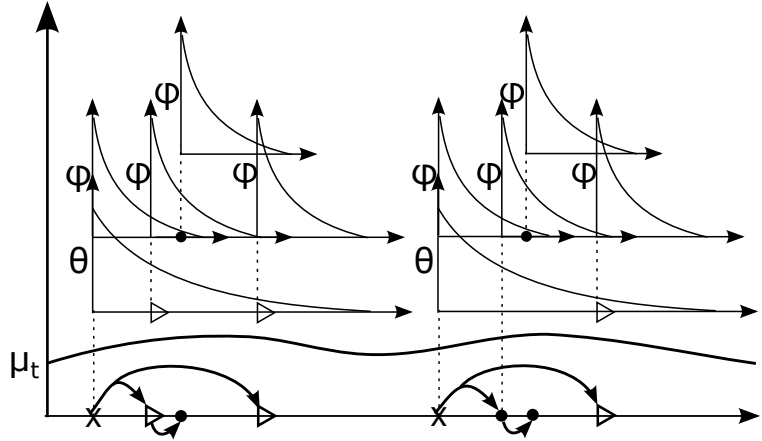}}
\caption{A realization of the ARMA point process with innovation (immigration) intensity $\mu$, MA (shot noise) intensity $\theta$, and AR (Hawkes) intensity $\phi$. Immigrants, $\theta$-offspring and $\phi$-offspring are denoted by x, triangle, and dot, respectively. A point is connected to the intensity that it triggers by a vertical dashed line. All points are projected onto the horizontal axis, with parenthood indicated by arrows, forming the full realization.}
\label{fig:Realisation}
\end{figure}
Let $(\Omega,\F,\Pm)$ be a probability space, let $\mu\in(0,\infty)$, let $N^\mu:\Omega\times\B(\R)\rightarrow\Z$ be a homogeneous Poisson process on $\R$ with rate $\mu$, let $\eta\in[0,1)$, $\gamma\in[0,\infty)$, let $\theta,\phi:[0,\infty)\rightarrow[0,\infty)$ be integrable functions with the property that 
\begin{equation}
\int_0^\infty \theta(t)dt=\gamma~~~ {\rm and}   ~~~\int_0^\infty \phi(t)dt=\eta~, 
\label{etjkukjrn}
\end{equation}
let $N^I:\Omega\times\R\rightarrow\Z$ be a cluster process generated by $N^\mu$ with i.i.d. clusters such that, given a cluster center $y\in\R$, a cluster $N_y$ is distributed according to  
\begin{enumerate}
\item a singular Dirac measure $\delta_y(\cdot)$ at $y$ (\emph{counting the immigrant}) and
\item additional points distributed according to an inhomogeneous Poisson process with rate $\theta(\cdot-y)$,
\end{enumerate}
then a point process $N:\Omega\times\B(\R)\rightarrow\Z$ is called \textit{ARMA point process} if it is the collection of Galton-Watson branching processes with Poisson offspring distribution $\phi$ based on immigrants $N^I$. That is, if both the original immigrants $N^\mu$ and their $\theta$-offspring trigger a Galton-Watson branching process with intensity $\phi$. This process is visualized in Figure~\ref{fig:Realisation}.

\subsection{Properties}
Lemma 6.3.II and Exercise 6.3.5 in \cite{VereJones2003_vol1} ensure the existence and stationarity of $N^I$ and $N$. The main point here is that, for a stationary immigration process, there exists a stationary cluster process if clusters are i.i.d. distributed (given the cluster center) and the cluster size is finite, which is ensured by $\eta\in[0,1)$ and $\gamma\in[0,\infty)$. The construction of a point process as a collection of Galton-Watson branching processes with Poisson offspring distribution was initially used in the point process representation of the Hawkes process in \cite{Hawkes1974}, see also Example 6.3(c) in \cite{VereJones2003_vol1}.

For the full filtration $\bm{\F}=\sigma(N^\mu,N)$, it holds that the $\bm{\F}$-conditional intensity function is given by
\begin{equation} \label{eq:armaP}
\lambda(t) 
= \mu + \int_{-\infty}^{t} \theta(t-s)dN^\mu_s + \int_{-\infty}^{t} \phi(t-s)dN_s. 
\end{equation}

With respect to the filtration $\bm{\F}=(\F_{t})_{t\in(-\infty,\infty)}$, the intensity function can be seen as a conditional hazard function in the sense that 
\begin{equation}
\lambda(t)=\text{lim}_{\Delta \downarrow 0}~\Delta^{-1}\text{E}\left[ N(t,t+\Delta) | \F_{t-} \right].
\end{equation}
For details see Chapter 7 of \cite{VereJones2003_vol1}. Note that the conditional intensity function uniquely defines the probability structure of a point process only if it is measurable with respect to its \emph{internal history}.\footnote{See Chapter 7 and Proposition 7.2.IV. in \cite{VereJones2003_vol1} for details.} Thus, unlike the case for the classical Hawkes process, the  $\bm{\F}$-conditional intensity function \eqref{eq:armaP} is not a defining property of $N$.  

Alternatively, one can use a random sequence of indicator variables depending on $N$ and $N^\mu$ to concisely describe the conditional intensity and the information flow. For $t\in\R$, we can define the random sequence $\mathcal{Z}_t=(Z_i)_{i\in\N}$ where, for $i\in\Z\cap(-\infty,N(t)]$, $Z_i=1$ if the $(N(t)-i+1)$-th last point was an immigrant and $Z_i=0$ otherwise. Then the full history is generated by $N$ and $\mathcal{Z}$, that is $\bm{\F}=\sigma(N^\mu,N)=\sigma(\mathcal{Z},N)$ and the $\bm{\F}$-conditional intensity can be written as 
\begin{as}\label{eq:armaPsum}
\lambda(t)
&= \mu + \sum_{j=-\infty}^{N^{\mu}(t)} \theta(t-T^{\mu}_j) + \sum_{k=-\infty}^{N(t)} \phi(t-T_k)~\\
&= \mu + \sum_{j=-\infty}^{N(t)} Z^\mu_j \theta(t-T_j) + \sum_{k=-\infty}^{N(t)} \phi(t-T_k)~,
\end{as}
where $(T^{\mu}_j)_{j\in\Z}$ and $(T_j)_{j\in\Z}$ denote the jump times of $N^{\mu}$ and $N$, respectively.

\subsection{First and second order statistics}
Using stationarity of the process, we can take the expectation of the $\bm{\F}$-conditional intensity function in \eqref{eq:armaP} to get the expected intensity,
\begin{equation}\label{eq:Eint}
\bar\lambda=\frac{ \mu(1+\gamma) }{ 1- \eta }~,
\end{equation}
which defines the \emph{first moment measure} of $N$, that is, for any Borel set $A\subseteq\R$, we have $E[N(A)]=\int_A \bar\lambda dt$.
By stationarity of the process, the \emph{covariance measure}, if it exists, is fully defined by its density $c(u)=Cov(dN_t,dN_{t-u})$,\footnote{Hereafter, we will use artificial objects like $Cov(dN_t,dN_{t-u})$ for simplicity. What we actually mean here is that $c$ is a density in the sense that for Borel sets $A,B\subseteq\R$ it holds that $Cov(N(A),N(B))=Cov(\int_{A}dN_t,\int_{B}dN_t)=\int_{A\times B}c(u-t)du dt$. For a rigorous definition of moment measures and their densities, we refer to Chapter 5.4 in \cite{VereJones2003_vol1}.} see 6.1.I. in \cite{VereJones2003_vol1}. The covariance measure has a singular Dirac component at $0$ and thus $c$ can be written as 
\begin{equation}\label{eq:covdens}
c(u)=\overline{\lambda}\delta(u)+\overline{\lambda}h(u)-\overline{\lambda}^2
\end{equation}
with the symmetric function $h:(-\infty,\infty)\rightarrow[0,\infty)$, called \textit{palm-intensity} in \cite{CoxIsham1980}, given by
\begin{as} 
\label{eq:PI}
h(u)&=\Pm\left[ dN_{t+u}=1\vert dN_t=1\right]\dfrac{1}{du}=\E\left[ dN_t dN_{t-u}\right]\dfrac{1}{\overline{\lambda}dt du},\quad u\in(0,\infty)\\
h(0)&=0.
\end{as}
 
To derive an expression for $h$ utilizing equation \eqref{eq:armaP}, we first need to derive an equation for the function $\tau:(-\infty,\infty)\rightarrow[0,\infty)$ given by the defining equation
\begin{equation}
\label{eq:tau1}
\E\left[ dN_t dN^{\mu}_{t-u}\right]=\left(\delta(u)\mu+\tau(u)\right)dt du.
\end{equation}
We multiply equation \eqref{eq:armaP} with ${dN^{\mu}_{t-u}}/{du}$ and take expectations to arrive at
\begin{as}
\tau(u)=\mu^2+\int_0^t \theta(t-s)\dfrac{1}{du}\E\left[ dN^{\mu}_s dN^{\mu}_{t-u}\right]+\int_{-\infty}^t\phi(t-s) \tau(s-t+u)ds,\quad u\neq 0.
\end{as}
Using $\E\left[ dN^{\mu}_s dN^{\mu}_{t-u}\right]=\left(\mu\delta(s-t+u)+\mu^2\right)ds du$ we get
\begin{as}\label{eq:tau2}
\tau(u)&=\mu^2(1+\gamma)+\mu(\theta(u)+\phi(u))+\int_0^\infty\phi(s)\tau(u-s)ds,\quad u\in(0,\infty),\\
\tau(0)&=\overline{\lambda}\mu.
\end{as}
Now we are concerned with calculating the full conditional intensity function $h$. To this end, let $\rho:[0,\infty)\rightarrow[0,\infty)$ be the function given by $\rho=\overline{\lambda}h$.
Multiplying equation \eqref{eq:armaP} with ${dN_{t-u}}/{du}$, taking expectations and using identities \eqref{eq:PI} and \eqref{eq:tau2} gives
\begin{as}
\label{eq:PIsol}
\rho(u)&=\mu\overline{\lambda}+\mu\theta(u)+\int_{-\infty}^{t-u}\theta(t-s)\tau(t-u-s)ds+\overline{\lambda}\phi(u)+\int_0^\infty\phi(s)\rho(u-s)ds\\
&=\mu\overline{\lambda}+\mu\theta(u)+\overline{\lambda}\phi(u)+\int_0^\infty\theta(s+u)\tau(s)ds+\int_0^u\theta(s)\tau(s-u)ds\\
&\quad+\int_0^\infty\phi(u+s)\rho(s)ds+\int_0^u\phi(s)\rho(u-s)ds.
\end{as}
If equation \eqref{eq:PIsol} has a solution, the covariance measure exists and has density $c$ defined in \eqref{eq:covdens}. As $\phi$ and $\theta$ are integrable functions, by $L^1$-theory for Fredholm integral equations, see, e.g., chapter 2.3 in \cite{Pipkin1991}, there exist locally integrable solutions to equations \eqref{eq:tau2} and \eqref{eq:PIsol}.
For exponential densities, equation \eqref{eq:PIsol} can be solved explicitly. To this end, let $\theta_0,\theta_1,\phi_0\in(0,\infty)$, $\phi_1\in(\phi_0,\infty)$ and let $\theta(t)=\theta_0 e^{-\theta_1 t}$ and $\phi(t)=\phi_0 e^{-\phi_1 t} = \eta \phi_1 e^{-\phi_1 t}$ (from the definition of the branching ratio $\eta$ given by (\ref{etjkukjrn}), which shows indeed that the condition $\eta < 1$ is equivalent to $\phi_1=\phi_0/\eta > \phi_0$). Then, as we show in appendix \ref{sec:app1}, there exist constants $K_1,K_2\in(0,\infty)$ depending on $(\theta_0,\theta_1,\phi_0,\phi_1)$ such that the palm intensity $h$ is given by
\begin{equation}\label{eq:renorm}
h(t)=\overline{\lambda}+K_1 e^{-(\phi_1-\phi_0) t}+K_2 e^{-\theta_1 t} = 
\overline{\lambda}+K_1 e^{-(1-\eta) \phi_1 t}+K_2 e^{-\theta_1 t},\quad t\in(0,\infty).
\end{equation}
The term $e^{-(1-\eta) \phi_1 t}$ recovers the standard renormalisation of the characteristic time scale ${1 \over \phi_1}$
to ${1 \over \phi_1}{1 \over 1-\eta}$  by counting over all generations of the $\phi$-process.

For the special cases of a Hawkes or a Neyman-Scott processes, the above expression reduces to the respective well-known palm intensities, in particular it holds that $K_1=0$ for $\phi_0=0$ and $K_2=0$ for $\theta_0=0$. As for the Hawkes and the NS processes, we have $\lim_{t\rightarrow\infty}h(t)=\overline{\lambda}$.

     
\subsection{Comments} \label{sec:comments} 
\paragraph{Branching interpretation.}
The explicit construction of the ARMA point process brings great insight, as visualised in Fig.~\ref{fig:Realisation}: $\mu$ introduces \emph{immigrants}, which may then trigger a single generation of $\theta$-\emph{offspring} with intensity $\theta(\cdot)$, and then all existing points ($\mu$-immigrants and their $\theta$-\emph{offspring}) trigger a generation of $\phi$-\emph{offspring} with intensity $\phi(\cdot)$, which may, in turn, trigger the subsequent generation of $\phi$-\emph{offspring} in the same way.
The sum of these independent inhomogeneous Poisson processes provides the $\bm{\F}$-conditional intensity \eqref{eq:armaP}, and the set of immigrant and ($\theta$- and $\phi$-)offspring points forms the ARMA point process realization. The \emph{branching ratio} $\gamma$ is the expected number of $\theta$-\emph{offspring} of a single immigrant, and the \emph{branching ratio} $\eta$ is the expected number of immediate $\phi$-\emph{offspring} of any point. Further, counting all generations, a single point is expected to produce $\eta+\eta^2+\dots=1/(1-\eta)-1$ $\phi$-\emph{offspring}. Thus, as in the Hawkes process, $\eta$ is the expected proportion of all points that are $\phi$-\emph{offspring}. From the \emph{ARMA point process} defined above, by setting $\gamma=0$ one recovers the Hawkes process, and by setting $\eta=0$ one recovers a modified Neyman-Scott (NS) process\footnote{I.e., an NS process where the immigrant is included in the counts.} (see \cite{Neyman1958}). It can thus be regarded as an extension of both a self-exciting and an i.i.d. cluster process. 

\paragraph{Extensions/Modifications.}
\begin{enumerate}
\item
The ARMA point process can be extended by using an inhomogeneous Poisson process $N^{\mu}$ in definition \ref{sec:def} for a locally integrable rate function $\mu\colon(-\infty,\infty)\rightarrow[0,\infty)$. The resulting cluster process of immigrants, $N^{I}$, can be extended to marked point process with the marks $(Y_i)_{i\in\N}$ being i.i.d. non-negative valued random variables, such that the cluster generated by the $j^{th}$ point of $N^\mu$ at $t_j\in\R$  is distributed according to an inhomogeneous Poisson process with rate $y_j\theta(\cdot-t_j)$ for a realisation $y_j$ of $Y_j$. Moreover, the ARMA point process $N$ can be extended to a marked point process with the marks $(Z_i)_{i\in\N}$ being i.i.d. non-negative valued random variables, such that the cluster generated by the $k^{th}$ point of $N$ at $t_k\in\R$  is distributed according to an inhomogeneous Poisson process with rate $z_k\theta(\cdot-t_k)$ for a realisation $z_k$ of $Z_k$.
Then the $\bm{\F}$-conditional intensity of $N$ in the form of \eqref{eq:armaPsum} can be written as
\begin{equation}\label{eq:armamarks}
\lambda(t)
= \mu(t) + \sum_{j=-\infty}^{N^{\mu}(t)} Y_j \theta(t-T^{\mu}_j) + \sum_{k=-\infty}^{N(t)} Z_k \phi(t-T_k)~.
\end{equation}
\item We will also show how to estimate a variant of the process \eqref{eq:armamarks} with background immigration, that is, where the immigrants $N^{\mu}$ are not included in the sample. In this case, the $\bm{\F}$-conditional intensity of the process takes the form 
(\ref{eq:armamarks}) without the first term $\mu(t)$ in the r.h.s.
\begin{equation}\label{eq:armaback}
\lambda(t)
=\sum_{j=-\infty}^{N^{\mu}(t)} Y_j \theta(t-T^{\mu}_j) + \sum_{k=-\infty}^{N(t)} Z_k\phi(t-T_k)~.
\end{equation}
While the inclusion of the immigrants can be argued to be \emph{natural} due to its similarity with the Hawkes process or INARMA time series (introduced below), which both count immigrants, their exclusion -- i.e., the use of background immigration -- can be argued to be \emph{natural} by interpreting it as a Hawkes process with Neyman-Scott immigration. In practice the application may clarify which specification makes sense.
\end{enumerate}

\paragraph{Hawkes processes with general immigrants in the literature.}
A framework for Hawkes processes with general immigrants that covers the ARMA point process has been introduced in \cite{Bremaud2002}, where they give an explicit form of the Bartlett \cite{Bartlett1963} spectral measure. 

One specific interpretation of the Hawkes process with general immigrants from \cite{Bremaud2002} has been introduced as the \emph{dynamic contagion process} in \cite{dassios2011dynamic}. The authors define a point process that has a stationary\footnote{In \cite{dassios2011dynamic}, they work in a non-stationary framework starting from an initial value at $\lambda_0$ at $0$, whose influence diminishes over time $t\rightarrow\infty$.} conditional intensity function, with respect to a suitable history, 
\begin{equation}\label{eq:DZ}
\lambda(t)
= \nu + \sum_{j=-\infty}^{N^{\mu}(t)} Y_j \theta(t-T^{\mu}_j) + \sum_{k=-\infty}^{N(t)} Z_k\phi(t-T_k)~.
\end{equation}
for $\nu\in[0,\infty)$ that corresponds to a (\say{background}) immigration Poisson process $N^{\nu}$, a Poisson process $N^{\mu}$ defined as above and independent from $N^{\nu}$ and (respectively) identically distributed marks $(Y_j)_{j\in\Z},(Z_k)_{k\in\Z}$. 
Thus, for the sake of comparison with \eqref{eq:armamarks}, we will set the background Poisson immigration rate $\nu=\mu$.
Observe that, in this special case, the conditional intensity function is the same as of the ARMA point process. 
However, the processes are distinct as there is a difference in the dependence structure between $N$ and $N^{\mu}$: for the ARMA point process, immigration and external shocks to the intensity are, as in the original Hawkes process, identical. The model \eqref{eq:DZ} from \cite{dassios2011dynamic} assumes an underlying immigration process with rate $\nu$ and an \emph{additional independent} process $N^\mu$ that leads to jumps of the intensity $\lambda$, but not of $N$.\footnote{By using a mixture distribution for the marks $(Y_i)_{i\in\N}$ in \eqref{eq:armamarks}, the ARMA point process can mimic the structure that some immigrants (cf. $N^{\nu}$) do not trigger a cluster but some do (cf. $N^{\mu}$). However, every immigrant ($N^{\nu}$ and $N^{\mu}$) will be counted. In this sense, a very similar but not the exact same process is covered in the present paper.} Let us also note that, for the special case $\nu=0$, the dynamic contagion process is identical to modification (2) of the ARMA point process mentioned above. 
In \cite{dassios2011dynamic}, it is assumed that the intensities $\theta$ and $\phi$ have an exponential form and decay at the same exponential rate $\delta$. This severe restriction leads to a very specific Markov structure and allows one to apply the powerful Markov process theory to calculate the Laplace transform of the process and its moments. 

By keeping track of the age of immigrants, the Markov property has recently been extended in \cite{Boumezoued2016} to a more general version of the dynamic contagion process, allowing for an explicit description of the Laplace transform for certain differentiable intensities (see assumption 1 in \cite{Boumezoued2016}). 

The next section provides a justification for our (different) choice of the specific structure of $N$ and the label \emph{ARMA point process}.


\section{The Relationship to Integer-Valued Time Series}\label{sec:inarma}
The ideas in this section are based on Kirchner \cite{kirchner2016hawkes}, who considers the framework of a classic Hawkes process.

\subsection{The INARMA model}
\begin{figure}[h!]
\centerline{\includegraphics[width=12cm]{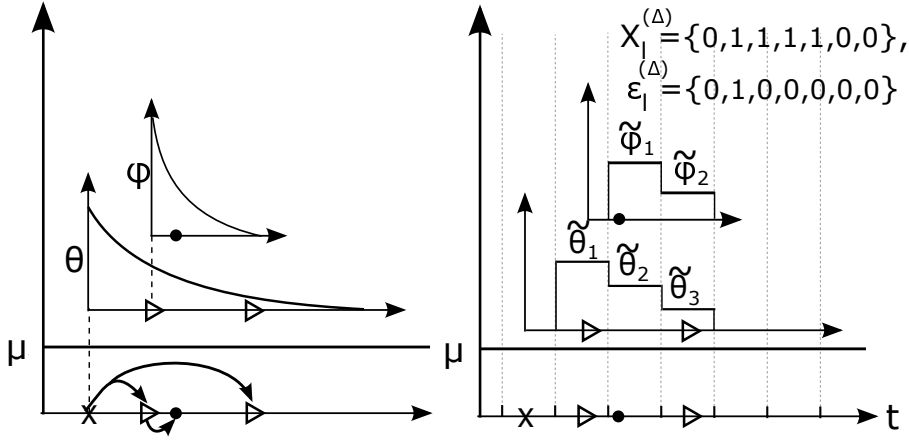}}
\caption{Representation of a cluster generated by the ARMA point process (left plot), and the INARMA$(2,3)$ process (right plot) which approximates the ARMA point process with grid $\Delta$. Only intensities that generated a point are shown. The immigrant and total counts are given which form the INARMA realization rather than the points. The origins of the axes framing the AR and MA triggering coefficients are located at the time values of the points that triggered them to highlight a source of approximation error: the INARMA can only trigger across bins, not within them. }
\label{fig:Analogues}
\end{figure}
The discrete valued analogues of classical time series models \cite{brockwell2013time} have seen a flurry of recent development \cite{Fokianos, McKenzie2003,weiss2008thinning,weiss2008serial} and enjoy many current and potential applications. Here we consider the INARMA$(p,q)$ process, an integer-valued ARMA process for counts. For its definition, we will need the \emph{Poisson thinning} operator $\circ$~ to preserve the count value of the process. In all of the following, we assume an underlying probability space $(\Omega,\F,\Pm)$. For some count variable $Z:\Omega\rightarrow\N$, a non-negative real number $\alpha\in[0,\infty)$, a sequence of i.i.d. Poisson distributed random variables $(Y_{i})_{i\in\Z}\overset{i.i.d}{\sim} \text{Pois}(\alpha)$, we define the count variable $\alpha \circ Z:\Omega\rightarrow\N$ by
\begin{equation}\label{eq:Thin}
\alpha \circ Z=\sum_{i=1}^{Z} Y_{i},~\alpha>0,~\text{and}~\alpha \circ 0 := 0.
\end{equation}

\begin{definition}\label{def:inarma}
Let $p,q\in\N$, let $\widetilde{\mu}\in[0,\infty)$, let $(\epsilon_i)_{i\in\Z}$ be a sequence of i.i.d. random variables $\epsilon_{i}\overset{i.i.d}{\sim} \text{Poisson}(\widetilde{\mu})$, let $(\widetilde{\theta}_k)_{k=1,\dots,p}$ and $(\widetilde{\phi}_j)_{j=1,\dots,q}$ be real-valued, non-negative sequences with the property that $\sum_{j=1}^{p}\widetilde{\phi}_j<1$. Then an integer-valued time series $X:\Omega\times\Z\rightarrow\N$ is an \emph{INARMA$(p,q)$ process} if it satisfies the difference equation
\begin{equation}\label{eq:INARMA}
X_{l}=\epsilon_{l} + \sum_{k=1}^{q}\widetilde{\theta}_{k} \circ \epsilon_{l-k}+ \sum_{j=1}^{p}\widetilde{\phi}_{j}\circ X_{l-j},  \quad l\in\Z,
\end{equation}
where all thinning operations are mutually independent.\footnote{There exist other interpretations for the serial dependence of thinning operations in the literature, see, e.g., chapter 5 of \cite{Fokianos02} for a specific interpretation in the general case or \cite{brannas2001estimation} for an overview of interpretations in the INARMA$(0,q)$-case.} 
\end{definition}

\subsection{Comments}\label{sec:com2}
The INARMA process \eqref{eq:INARMA} exists as a multi-type branching process that is stationary for $\sum_{j=1}^{p}\widetilde{\phi}_j<1$, see corollary 2 in \cite{Latour95}. 

The autocovariance function $\gamma:\Z\rightarrow[0,\infty)$ of the INARMA process \eqref{eq:INARMA} is the function $\gamma(u)=Cov(X_l,X_{l-u})$, often equivalently used via the rescaled autocorrelation function (ACF). It is explicitly understood for INAR$(p)$ and INMA$(q)$ processes with Bernoulli thinning, see \cite{DuLi91} and \cite{McKenzie1988}, respectively. The case of Poisson thinning is covered by a general result in \cite{Latour98} for INAR$(p)$ and can be extended directly from \cite{McKenzie1988} for INMA$(q)$ processes. To our knowledge, an explicit description of the ACF for the full INARMA$(p,q)$ process \eqref{eq:INARMA} has not been established yet and is left for future research.
Given the autocovariance $\gamma$, we can find best linear conditional predictors for $X$ in the space of real-valued time series and one might extend the notion of partial autocorrelation (PACF) to INARMA$(p,q)$-processes.\footnote{See, e.g., sections 3.4 and 5.2 \cite{brockwell2013time} for a derivation of the partial autocorrelation for real-valued time series in terms of optimal linear predictors, using the Durbin-Levinson algorithm.} 

From the difference equation \eqref{eq:INARMA}, assuming a time distance $\Delta\in[0,\infty)$ between the counts, it is straightforward to deduce a discrete conditional intensity function $\lambda^{(\Delta)}$ given by
\begin{equation}\label{eq:INARMAint}
\lambda^{(\Delta)}(l)=\dfrac{1}{\Delta}E[X_l|X_{(l-1):(l-p)},\epsilon_{(l-1):(l-q)}]=\dfrac{1}{\Delta}\left( \widetilde{\mu} + \sum_{k=1}^{q}\widetilde{\theta}_{k} \epsilon_{l-k}+ \sum_{j=1}^{p}\widetilde{\phi}_{j} X_{l-j}\right),\quad l\in\Z.
\end{equation}

The branching interpretation of the INARMA process \eqref{eq:INARMA} is the same as for the ARMA point process: the innovation count $\epsilon_l$ introduces immigrants, and thinning~\eqref{eq:Thin} has the interpretation that each of the $X_{j}$ points in the $j^\text{th}$ bin is expected to produce $\widetilde{\phi}_{l-j}$ offspring in the $l^\text{th}$ bin, where $l>j$. Thus the INARMA process introduces a burst of offspring triggered by immigrants, via thinning with $\widetilde{\theta}$ coefficents, and an autoregressive tree of offspring triggered by all past events, via thinning with $\widetilde{\phi}$ coefficients. 

Observe that the thinning operator defined in equation \eqref{eq:Thin} has the property that $\alpha \circ Z|Z=z$ is a sum of $z$ independent Poisson variables with parameter $\alpha$, and has distribution $\text{Pois}(\alpha z)$. Thus, given that both all thinnings in \eqref{eq:INARMA} are independent of each other, and of the Poisson innovation $\epsilon_l$, then the conditional df of $X_l|X_{(l-p):(l-1)},\epsilon_{(l-q):(l-1)}$ is also Poisson. The unconditional df of $X_l$, on the other hand, is not Poisson. In this sense, we maintain structural similarities to ARMA point processes, whose distributions, given the full history, are that of an inhomogeneous Poisson process. It is important to note that the standard thinning used in integer time series is Binomial thinning, where the variable $Y$ has a Bernoulli distribution. In this case, the unconditional distribution of $X_l$ is Poisson, but the conditional one is not. A survey of the different thinning specifications employed within the literature are summarized in \cite{weiss2008thinning}. 


\subsection{Connection to the ARMA point process}    
As discussed above, the branching interpretation of the INARMA time series and the ARMA point process are identical. Next, a formal argument is made to suggest an asymptotic equivalence between the process from section \ref{sec:def} and INARMA time series, motivating the term \emph{ARMA point process}. For this, let $N$ be an ARMA point process as defined in section \ref{sec:def} given by an immigration rate $\mu\in(0,\infty)$ and integrable intensities $\theta,\phi:[0,\infty)\rightarrow[0,\infty)$. 
If one aggregates the ARMA point process $N$ on bins of width $\Delta>0$, one obtains the counting variables $\lbrace N^{(\Delta)}_{l}=N\left(\Delta l,\Delta (l+1) \right),~l\in\Z\rbrace$ for all points and  $\lbrace \epsilon^{(\Delta)}_{l}=N^\mu\left(\Delta l,\Delta (l+1) \right),~l\in\Z\rbrace$ for innovations, see Fig.~\ref{fig:Analogues} for an example. 
Then, for mild assumptions on $\theta$ and $\phi$, we $\Pm$-a.s. have the convergence 
\begin{as}\label{eq:pwconv}
\sum_{k=1}^{q} \theta(k\Delta)\epsilon^{(\Delta)}_{\floor{t/\Delta}-k}&\underset{p,q\rightarrow\infty,\Delta\rightarrow 0}{\longrightarrow}\int_{-\infty}^{t} \theta(t-s) dN_s^\mu,\\
\sum_{j=1}^{p} \phi(j\Delta ) N^{(\Delta)}_{\floor{t/\Delta}-j}&\underset{p,q\rightarrow\infty,\Delta\rightarrow 0}{\longrightarrow}\int_{-\infty}^t \phi(t-s) dN_s.
\end{as}
Thus, for the aggregated model $N^{(\Delta)}$, using \eqref{eq:armaP} and \eqref{eq:pwconv}, we expect a discrete conditional intensity function
\begin{as}\label{eq:approx}
\lambda^{(\Delta)}(l)=
\mu + \sum_{k=1}^{q}\theta(k \Delta ) \epsilon^{(\Delta)}_{l-k} + \sum_{j=1}^{p} \phi(j\Delta ) N^{(\Delta)}_{l-j}+\err(\Delta,p,q)~,
\end{as}
where  $\err(\Delta,p,q)$ consists of approximation of the integral in \eqref{eq:pwconv} with finite step functions of length $p,q$ and the influence of offsprings that are triggered by points $N$ or $N^\mu$ in the \emph{same} interval (of length $\Delta$). Thus, we expect $\err(\Delta,p,q)\rightarrow 0$ for $p,q\rightarrow\infty$ and $\Delta\rightarrow 0$ and the aggregated ARMA point process to (approximately) follow an INARMA process. This formal reasoning leads us to the conjecture that the finite dimensional distributions of a sequence of INARMA processes converge to the finite dimensional distributions of the ARMA point process. 

\begin{conjecture}
\begin{enumerate}
\item \emph{(The INARMA$(\infty,\infty)$ process.)}
Let $\widetilde{\mu}\in[0,\infty)$, let $(\epsilon_i)_{i\in\Z}$ be sequence of i.i.d. random variables $\epsilon_{i}\overset{i.i.d}{\sim} \text{Poisson}(\widetilde{\mu})$ and let $(\widetilde{\theta}_k)_{k\in\N}$ and $(\widetilde{\phi}_j)_{j\in\N}$ be real-valued, non-negative sequences with the property that $\sum_{k=1}^{\infty}\widetilde{\theta}_k<\infty$ and $\sum_{j=1}^{\infty}\widetilde{\phi}_j<1$. Then there exists an integer-valued stationary time series $X:\Omega\times\Z\rightarrow\N$ that satisfies the difference equation
\begin{equation}\label{eq:INARMAinf}
X_{l}=\epsilon_{l} + \sum_{k=1}^{\infty}\widetilde{\theta}_{k} \circ \epsilon_{l-k}+ \sum_{j=1}^{\infty}\widetilde{\phi}_{j}\circ X_{l-j},  \quad l\in\Z.
\end{equation}
\item \emph{(Approximation of the ARMA point process.)} Let $N$ be an ARMA point process as defined in section \ref{sec:def} with a piecewise continuous intensity $\phi$, then there exists $\delta\in(0,\infty)$ such that 
\begin{itemize}
\item for $\Delta\in(0,\delta)$ equation \eqref{eq:INARMAinf} with $\widetilde{\mu}=\mu\Delta$, $\widetilde{\theta_k}=\Delta\theta(k\Delta)$, and $\widetilde{\phi_j}=\Delta\phi(j\Delta)$ defines a stationary INARMA process $X^{(\Delta)}$ and
\item the family of point processes $(N^\Delta)_{\Delta\in(0,\delta)}$ given by 
\begin{equation}
N^\Delta(A)=\sum_{n\colon n\Delta\in A}X_n^\Delta,\quad \text{for a Borel set }A\subseteq\R
\end{equation}
converges weakly\footnote{\emph{Weak convergence} of point processes is understood as vague convergence of their induced measures and equivalent to a convergence of their finite dimensional distributions, see, e.g., \cite{Resnick1987}.} to $N$ for $\Delta\rightarrow 0$.
\end{itemize}
\end{enumerate}
\end{conjecture}   

A rigorous proof goes beyond the scope of this article. As special cases, the aggregated Hawkes process is approximated by the INAR process, and the aggregated Neyman-Scott process is approximated by an INMA process. The weak convergence of the INAR process to the Hawkes process was established in \cite{kirchner2016hawkes}. 
 
\subsection{Implications}  
In case the above conjecture, as suggested by the formal argument, turns out to be true, it verifies the already useful analogy between time series models and (binned) point processes. 

As an important example, the autocorrelation (ACF) and partial autocorrelation (PACF) functions, which have been thoroughly studied for real-valued time series and are extendible to integer-valued time series (see the comments in section \ref{sec:com2} above), can be used to characterize ARMA point processes. The ACF, being defined as the covariance between two lagged random variables for time series, is directly related to the palm intensity \eqref{eq:PI}, which defines the covariance measure of the point process. An analogy to the PACF, whose definition involves the notion of linear (best) predictors given a sub-$\sigma$-algebra, is not as readily defined for a point process, and leaves room for further research.

Moreover, as is the case for the ARMA point process, the INARMA model cannot be directly/ simultaneously estimated by MLE due to missing information -- here being the innovation counts $\epsilon^{(\Delta)}_{l-1},\epsilon^{(\Delta)}_{l-2},\dots$ where only the complete counts $X^{(\Delta)}_{l-1},X^{(\Delta)}_{l-2},\dots$ are observed. The EM algorithm provided in section \ref{sec:EM} to estimate the ARMA point process may also be applied to the INARMA model and thus provides a powerful approach to fitting INARMA time series models.  


\section{Estimation of the ARMA point process with EM algorithms}\label{sec:EM}
\subsection{Motivations for the EM scheme}
 Unlike the Hawkes process, the conditional intensity functions of the NS, the ARMA, and marked extensions \eqref{eq:armamarks} depend on (knowing) the immigrants. These processes are Poisson, given the conditional intensity, which is itself stochastic. For such processes, in general \cite{moller2003statistical}, the likelihood is given in terms of an expectation with respect to the unobserved random intensity function. The typical solution in this case is to perform likelihood inference by MCMC (Markov Chain Monte Carlo) sampling \cite{moller2003statistical}. In the case of univariate temporal point processes, this turns out to be rather simple to implement, and works well.
  
  Also, practically speaking, it is crucial in applications to consistently estimate trends in immigration $\mu(t)$ to avoid mistaking deterministic trends for stochastic fluctuations. For reasonable models with such features, moment-based estimation is not useful, and ``direct'' MLE via numerical maximization of the loglikelihood may perform poorly due to joint estimation of a large number of parameters perhaps along with non-parametric $\mu(t)$. In such a setting, the EM algorithm \cite{mclachlan2007algorithm} for maximum likelihood estimation is powerful as it decomposes otherwise unwieldy multi-parameter estimation into sub-estimations -- specifically here into very simple problem of density estimations from iid samples with weights. 
    
  Regarding scope, this EM framework for ARMA point-processes allows for the estimation of a range of model specifications, e.g.:
  \begin{itemize} 
   \item Submodels, including the Hawkes, NS, and SNCP.
   \item Immigrant points are observed (included in the sample) or unobserved, cf. \eqref{eq:armaback}.
   \item Trends in immigration and/or the branching ratios (as well as other parameters).
   \item Without or with marks (as in \eqref{eq:armamarks}) that are either observed, or unobserved and thus must be simulated within the EM algorithm.
   \item And, last but not least, INARMA time series with Poisson thinning \eqref{eq:INARMA} can be estimated by a simple and obvious modification of the EM algorithm. 
  \end{itemize}
  
  
Below, the EM algorithm for the case of the ARMA with marks and inhomogeneous immigration will be presented. Section \ref{sec:immknown} will be concerned with the simpler case where immigrants and marks are observed, section \ref{sec:immunknown} will treat the case of unobserved immigrants. Section \ref{sec:em_stepwise} concludes with a step-wise decription of the algorithm and section \ref{sec:em_comments} discusses convergence properties.

\subsection{Derivation of the the EM algorithm}\label{sec:deriv}
\subsubsection{Notations}
Let $T$ be the ARMA point process introduced in section \ref{sec:def} with observed realization $T_t(\omega)=\boldsymbol{t}=(t_1,\dots,t_n)$ on a fixed time window $[0,t]$. We will allow for inhomogeneous immigration with intensity $\mu(\cdot)$ and marks $(Y_i)_{i\in\N}$, cf. equation (\ref{eq:armamarks}). Then $T$ has $\sigma(N^\mu,N)$-conditional intensity
\begin{equation}
\lambda(s)
= \mu(s) + \sum_{j=-\infty}^{N^{\mu}(s)} Y_j \theta(s-T^{\mu}_j) + \sum_{k=-\infty}^{N(s)} \phi(s-T_k)~.
\end{equation}
For ease of presentation, we will introduce some notation.
\begin{enumerate}
\item
Let $C_t$ be the random sequence induced by the marked immigrant process $N^\mu$, with realization $\boldsymbol{c}=(c_1,\dots,c_{n_c})$ of length $N^{\mu}_t=n_c$. Each element has the time and mark $c_j=(s_j,y_j)\in(0,t]\times [0,\infty)$. Separately, denote the points  $\boldsymbol{s}=(s_1,\dots,s_{n_c})$ that are iid on window $(0,t]$ with density $\mu(\cdot)/\mu((0,t])$, and the marks as $\boldsymbol{y}=(y_1,\dots,y_{n_c})$, where marks are iid with pdf $m(\cdot)$. 
\item Similarly, denote by $\boldsymbol{o}=(o_1,\dots,o_{n_o})$ the realisation of the offspring process with $n_o=n-n_c$ and $\boldsymbol{t}=\boldsymbol{c}\cup\boldsymbol{o}$.
\item
Finally, let $f,g:[0,\infty)\rightarrow[0,\infty)$ be the density function of the AR and MA kernels, that is, $\theta=\gamma g$ and $\phi=\eta f$. Due to the presence of marks, without loss of generality we can assume that $\gamma\equiv 1$.
\end{enumerate}

\subsubsection{Description of the EM algorithm for observed immigrants and marks}\label{sec:immknown}
We will start with the simpler but less realistic case of where we observe a realisation of the random sequence $C_t(\omega)=\boldsymbol{c}=(c_1,\dots,c_{n_c})$ and show how the EM algorithm can be applied in this situation to derive parameter estimates based on a conditional expectation of the full likelihood (E-step) for the parameter vector $$\boldsymbol{\beta}=(\mu,f,m,g,\eta).$$ 
The parametric form or non-parametric subclass of the densities $\mu,f,m,g$ has to be specified in the maximisation step (M-step) to yield a well-defined maxisation problem. 

\paragraph{Likelihood expectation (E-step)}
To derive a likelihood function, we will use the fact that, given the full branching structure $\boldsymbol{Z}$, the ARMA process decomposes into three independent inhomogeneous Poisson processes (density $p_{\rho}$ for intensity $\rho$) and the density of the marks $m$. In particular, the density factorizes
\begin{equation}\label{eq:density}
p(\boldsymbol{t}\vert \boldsymbol{Z})=p_\mu(\boldsymbol{s})p_m(\boldsymbol{y})p_{\theta}(\boldsymbol{o}_{\theta})p_{\phi}(\boldsymbol{o}_{\phi}),
\end{equation}
where $\boldsymbol{o}_{\phi}$ and $\boldsymbol{o}_{\theta}$ denote the AR and MA offsprings, respectively. The branching structure $\boldsymbol{Z}$ is therefore a highly useful unknown, and will be treated as our EM ``missing data''. Formally define the missing data by indicator variables
\begin{equation}\label{eq:Z}
\boldsymbol{Z}=\lbrace {Z}^{\theta}_{i,j}~,~i=1,...,n,~j=1,...,n_c \rbrace \cup \lbrace {Z}^{\phi}_{i,j}~,~i=1,...,n,~j=1,...,n \rbrace~,
\end{equation}
which are zero except $Z^{\theta}_{i,j}=1$ if $t_i$ is triggered by $\theta(\cdot-s_j)$, and $Z^{\phi}_{i,j}=1$ if $t_i$ is triggered by $\phi(\cdot-t_j)$.

Given the missing data \eqref{eq:Z}, the ``complete data'' likelihood can be derived by using the missing data variables to rewrite \eqref{eq:density},
\begin{as}
    &L(\boldsymbol{\beta}\mid \boldsymbol{t},\boldsymbol{c},\boldsymbol{Z}) = \prod_{i=1}^{n_c} \mu(s_i) \text{Exp}\lbrace - \int_0^t \mu(s)ds \rbrace m(y_i) \times  \\
    & ~~~\prod_{i=1}^{n}\prod_{j=1}^{n_c} \Big[ y_j \theta(t_i-s_j) \Big]  ^{ Z_{i,j}^\theta } \text{Exp}\lbrace -   \sum_{j=1}^{n_c} y_j \int_0^t \theta(s-s_j) ds \rbrace \times \\
    & ~~~~~~\prod_{i=1}^{n}\prod_{k=1}^{n} \Big[ \phi(t_i-t_k) \Big] ^{Z_{i,k}^\phi} \text{Exp}\lbrace - \sum_{j=1}^{n}\int_0^t \phi(s-t_j)ds \rbrace~,\label{eq:likIntPart}
\end{as}
where an intensity is only evaluated at the times agreeing with the branching structure encoded within the missing data. 
Instead of optimizing the (inaccessible) complete data likelihood \eqref{eq:likIntPart}, the EM algorithm  uses the objective function
\begin{equation} \label{eq:Q} 
Q(\boldsymbol{\beta}\mid\boldsymbol{\widehat{\beta}})=\mathbb{E}_{\boldsymbol{Z}|\boldsymbol{\widehat{\beta}}, \boldsymbol{t},\boldsymbol{c} }
[ \log{L(\boldsymbol{\beta}\mid \boldsymbol{t},\boldsymbol{c},\boldsymbol{Z})} ]~,
\end{equation}    
the expectation of the log-likelihood over the missing data, given the complete observations $(\boldsymbol{t},\boldsymbol{c})$, and a parameter estimate $\boldsymbol{\widehat{\beta}}$, where at the $r+1^{th}$ iteration,
\begin{equation}\label{eq:EM-step}
\widehat{\boldsymbol{\beta}}^{(r+1)}=\max_{\boldsymbol{\beta}}Q(\boldsymbol{\beta}\mid\boldsymbol{\widehat{\beta}^{(r)}})~.
\end{equation}
The function Q then contains the probabilities
\begin{as}\label{eq:piPhi}
\pi^{\theta}_{i,j}&=\Pm\lbrace Z_{i,j}^\theta=1 |\boldsymbol{t},\boldsymbol{c},\boldsymbol{\beta}\rbrace=
\begin{cases}
\dfrac{ y_j\theta(t_i-s_j) }{ \sum_{j=1}^{n_c}  y_j\theta(t_i-s_j)  + \sum_{j=1}^{n}\phi(t_i-t_j)  }& \text{for } t_i\in\boldsymbol{o}\\
0 & \text{else}
\end{cases}\\
\pi^{\phi}_{i,j}&=\Pm\lbrace Z_{i,j}^\phi=1|\boldsymbol{t},\boldsymbol{c},\boldsymbol{\beta}\rbrace=
\begin{cases}
\dfrac{\phi(t_i-t_j)}{\sum_{j=1}^{n_c} y_j \theta(t_i-s_j) + \sum_{j=1}^{n}\phi(t_i-t_j)} & \text{for } t_i\in\boldsymbol{o}\\
0 & \text{else}
\end{cases}
\end{as}
in place of the missing data indicator variables \eqref{eq:Z}, while decoupling of the components in the likelihood \eqref{eq:likIntPart} is preserved. These probabilities follow from the thinning~\cite{VereJones2003_vol1} whereby the probability that $t_i$ comes from one of the independent (sub-)processes is equal to that process' share of the total conditional intensity function at $t_i$. For specified model components with given parameter estimate $\boldsymbol{\widehat{\beta}}$, these probabilities \eqref{eq:piPhi} can be computed.
    
\paragraph{Maximisation (M-step)}    
The maximisation of the expected log-likelihood \eqref{eq:Q} has the structure of probability density estimation with sample weights \eqref{eq:piPhi} in place of the indicator variables encoding the missing data. This decoupling into iid density estimation enables the estimation of relatively complex, as well as non-parametric densities. Specifically, 
\begin{enumerate}
\item $g$ is estimated on iid positive interevent times $\lbrace~ t_i-s_j,~ s_j<t_i ~\rbrace$ with weights $\pi^{\theta}_{i,j}$,
\item $f$ is estimated on iid positive interevent times $\lbrace~ t_i-t_j, ~t_j<t_i~ \rbrace$ with weights $\pi^{\phi}_{i,j}$,
\item $\mu/\mu((0,t])$ is estimated on immigration times $\boldsymbol{s}$,\footnote{To recover the estimated intensity, this density is multiplied by $n_c$ to satisfy $\int_0^{n_c}\mu(s)ds=n_c$.}
\item $m$ is estimated on the iid sample $\boldsymbol{y}$.
\end{enumerate}
The branching ratio estimator $\widehat{\eta}$ is
\begin{equation}\label{eq:etaEMarma}
	  \widehat{\eta}=\frac{ \sum_{i,j} \pi^\phi_{i,j} }{ \sum_{j=1}^{n} \int_0^{t-t_{j}}f(s)ds  }~,
\end{equation}
with cumulative distribution in the denominator to correct for expected offspring truncated by the end of the observation window.
    
The E-step and M-step may then be iterated, and the estimates taken when the parameter estimates and log-likelihood have converged.

\subsubsection{MCMC extension for unknown immigrants and marks}\label{sec:immunknown}
If we want to apply the EM algorithm in the more realistic case where we do not observe a sample $\boldsymbol{c}$ of the immigrant process $C$ -- i.e., we neither know which points are immigrants nor the values of their marks -- the missing data \eqref{eq:Z} becomes $(\boldsymbol{Z},C)$ and the objective function \eqref{eq:Q} takes the form
\begin{as} \label{eq:QM} 
Q(\boldsymbol{\beta}\mid\boldsymbol{\widehat{\beta}})&=
\mathbb{E}_{\boldsymbol{Z},C|\boldsymbol{\widehat{\beta}}, \boldsymbol{t} }
[ \log{L(\boldsymbol{\beta}\mid \boldsymbol{t},C,\boldsymbol{Z})} ]~\\
&=\mathbb{E}_{C |\boldsymbol{\widehat{\beta}}, \boldsymbol{t} }
\left[\mathbb{E}_{\boldsymbol{Z}|\boldsymbol{\widehat{\beta}}, \boldsymbol{t},C }
[ \log{L(\boldsymbol{\beta}\mid \boldsymbol{t},C,\boldsymbol{Z})} ]\right].
\end{as}    
Assume we are able to generate samples $\boldsymbol{c}^{(1)},\dots,\boldsymbol{c}^{(K)}$ from a density $p_C(\cdot\mid\boldsymbol{\widehat{\beta}}, \boldsymbol{t})$, we can then approximate the outer expectation with respect to $C$ to approximate the objective function by
\begin{as}\label{eq:obj_fct}
Q(\boldsymbol{\beta}\mid\boldsymbol{\widehat{\beta}})&\approx
\dfrac{1}{K}\sum_{k=1}^{K}\left[\mathbb{E}_{\boldsymbol{Z}|\boldsymbol{\widehat{\beta}}, \boldsymbol{t},\boldsymbol{c}^{(k)} }
[ \log{L(\boldsymbol{\beta}\mid \boldsymbol{t},\boldsymbol{c}^{(k)},\boldsymbol{Z})} ]\right].
\end{as}
To evaluate the inner expectation, we assume immigrants and marks are known ($C=\boldsymbol{c}^{(k)}$) and the procedure of section \ref{sec:immknown} applies.
In particular, for the probabilities $\pi^{\phi,k}_{i,j}$ and $\pi^{\theta,k}_{i,j}$ calculated according to \eqref{eq:piPhi} for $c^{(k)}$, then 
\begin{as}
&Q(\boldsymbol{\beta}\mid\boldsymbol{\widehat{\beta}}) \approx\frac{1}{K}\sum_{k=1}^K \sum_{i=1}^{n_c(k)} \log\left(\mu\left(s^{(k)}_i\right)\right) - \int_0^t \mu(s)ds + \frac{1}{K}\sum_{k=1}^K \sum_{i=1}^{n_c(k)} \log\left(m\left(y^{(k)}_i\right)\right) +  \\
    & ~~~\frac{1}{K}\sum_{k=1}^K\sum_{i=1}^{n}\sum_{j=1}^{n_c(k)} \pi^{\theta,k}_{i,j}\log\Big[ y^{(k)}_j \theta\left(t_i-s^{(k)}_j\right) \Big]  -   \frac{1}{K}\sum_{k=1}^K\sum_{j=1}^{n_c(k)} y^{(k)}_j \int_0^t \theta\left(s-s_j^{(k)}\right) ds\,  + \\
    & ~~~\sum_{i=1}^{n}\sum_{j=1}^{n}\left(\frac{1}{K}\sum_{k=1}^K\pi^{\phi,k}_{i,j}\right)\log\Big[ \phi(t_i-t_j) \Big] - \sum_{j=1}^{n}\int_0^t \phi(s-t_j)ds ~.\label{eq:obj_fct_detail}
    \end{as}
Thus, as before in section \ref{sec:immknown}, maximisation over parameters  $\boldsymbol{\beta}$ is decoupled into (weighted) iid density estimation and $g,f,\mu,m$ and $\eta$ can be estimated separately. The estimation of $\theta$ involves a pooling of interevent times $\{t_i-s_j^{(k)}\}$ for $k=1,\dots,K$; while the estimation of $\phi$ merely requires weights averaged over the ensemble.. 
However, the density of the immigrants
\begin{equation}\label{eq:dens}
p(\boldsymbol{c}\mid\boldsymbol{t}) = 
\frac{p(\boldsymbol{t}\mid \boldsymbol{c})p(\boldsymbol{c})}{p(\boldsymbol{t})}
\end{equation}
is not known analytically due to the lack of an expression for the denominator of \eqref{eq:dens}. 

\paragraph{Conditional simulation of MCMC sample.}
Instead, we employ a simple MCMC algorithm to simulate immigrant realizations from this distribution. The algorithm is extended from \cite{Moller2003} and \cite{moller2003statistical} (sections 7.1.2 and 10.2.1) and uses a Metropolis Hastings algorithm to generate a Markov chain sample from the unnormalized density $p_{C\mid T}(\cdot\mid \boldsymbol{t})p_T(\boldsymbol{t})$. 

Note that the joint probability $p(\boldsymbol{t},\boldsymbol{c})=0$ if $\boldsymbol{c}\not\subseteq\boldsymbol{t}$ and thus
\begin{equation}\label{eq:imm_density_cond}
p(\boldsymbol{c}\mid\boldsymbol{t})=K_1(\boldsymbol{t}) p\left(\boldsymbol{t}\mid\boldsymbol{c}\right)\prod_{s_i\in\boldsymbol{s}}\mu(s_i)=K_2(\boldsymbol{t}) p\left(\boldsymbol{o}\mid\boldsymbol{c}\right)\prod_{s_i\in\boldsymbol{s}}\mu(s_i),
\end{equation}
with constants $K_1,K_2$ depending on $\boldsymbol{t}$. Given the marked immigrants, the offspring $\boldsymbol{o}$ are a realization from a modified ARMA (cf. section \ref{sec:comments}) with (realized) conditional intensity
\begin{equation}
\lambda_{O\mid C}(s)
\approx \sum_{s_j\in\boldsymbol{s}} y_j \theta(s-s_j) + \sum_{t_j\in\boldsymbol{t}} \phi(s-t_j)~,
\end{equation}
ignoring the influence of points on $(-\infty,0)$.
Then, using proposition 7.2 III in \cite{VereJones2003_vol1},
\begin{equation}\label{eq:off_density}
p\left(\boldsymbol{o}\vert\boldsymbol{c}\right)=\prod_{t_i\in\boldsymbol{o}}\lambda_{O\mid C}(t_i)\text{Exp}\lbrace -\int_{0}^{t} \lambda_{O\mid C}(s)ds\rbrace.
\end{equation}	

\paragraph{Metropolis-Hastings iteration.}
Given the current state of the Markov chain, $\boldsymbol{c}^{(k)}=\boldsymbol{c}=(\boldsymbol{s},\boldsymbol{y})$, our Metropolis-Hastings iteration consists of proposing either the birth or death of an immigrant with probability $1/2$. 
In the case of birth, choose $s^\ast\in\boldsymbol{o}$ uniformly among offsprings $\boldsymbol{o}=\boldsymbol{t}\setminus\boldsymbol{s}$ and generate a mark $y^\ast$ from $m(\cdot)$.
Then, using equations \eqref{eq:imm_density_cond}-\eqref{eq:off_density} and $c^*=(s^*,y^*)$, the Metropolis-Hastings birth ratio\footnote{C.f. equation (7.6) in \cite{ moller2003statistical}.} reads
\begin{as}\label{eq:birth_ratio}
r_b(\boldsymbol{c},c^{*})&=\dfrac{n-n_c}{n_c+1}\dfrac{p(\boldsymbol{c}\cup c^*\mid t)}{p(\boldsymbol{c}\mid t)}\\
&=\frac{n-n_c}{n_c+1} \mu(s^{*})\text{Exp}\lbrace -\int_{0}^{t} y^{*}\theta(s-s^{*})ds \rbrace \times\\&
\quad\prod_{\substack{t_i\in \boldsymbol{o}\\t_i\neq s^*}} \left(  1+ \frac{y^{*}\theta(t_i-s^{*}) }{ \sum_{s_j\in \boldsymbol{s}} y_j\theta(t_i-s_j) + \sum_{t_j\in \boldsymbol{t} }\phi(t_i-t_j) } \right) \Bigg/  
\left( \sum_{s_j\in \boldsymbol{s} }y_j\theta(s^{*}-s_j) + \sum_{t_j\in \boldsymbol{t} }\phi(s^*-t_j)  \right),
\end{as}
and we accept and move to the new state $\boldsymbol{c}^{(k+1)}=\boldsymbol{c}\cup c^*$ with \emph{acceptance probability} $\min\{1,r_b(\boldsymbol{c},c^{*})\}$.
In the case of death, we choose $s^\ast\in\boldsymbol{s}$ uniformly and arrive at a death ratio
\begin{equation}\label{eq:death_ratio}
r_d(\boldsymbol{c},c^{*})=\dfrac{1}{r_b(\boldsymbol{c}\setminus c^{*},c^*)}
\end{equation}
and move to the new state $\boldsymbol{c}^{(k+1)}=\boldsymbol{c}\setminus c^*$ with \emph{acceptance probability} $\min\{1,r_d(\boldsymbol{c},c^{*})\}$.

\newpage
\subsection{Step-wise description of the EM algorithm}\label{sec:em_stepwise}
We use the notation of section \ref{sec:deriv}.
\vskip 0.5cm
 \setlength{\parindent}{0cm}
\begin{mdframed}[leftmargin=0pt,rightmargin=0pt]
0. Start with initial parameter "guess" $\boldsymbol{\beta}^{(0)}$ and a set of immigrants $\boldsymbol{c}^{(0)}\subseteq\boldsymbol{t}$.\\
\textit{For every iteration $r\in\N$, repeat E and M step as follows.}\\
\textbf{I. E-Step (MCMC)} \\
\textit{Generate Markov chain of immigrants $\boldsymbol{c}^{(r,1)},\dots,\boldsymbol{c}^{(r,K)}$ given $\boldsymbol{\beta}^{(r-1)}$. In every iteration $k$}\\
1. Flip a coin to choose birth or death.\\
2a. In case of birth, choose $\boldsymbol{s^*}\in\boldsymbol{o}^{(r,k-1)}$ uniformly and generate $y^\ast$ from $m^{(r-1)}$, calculate the birth ratio $r_b(\boldsymbol{c}^{(r,k-1)},c^{*})$ in \eqref{eq:birth_ratio} and accept $\boldsymbol{c}^{(r,k)}=\boldsymbol{c}^{(r,k-1)}\cup(s^*,y^*)$ with probability $\min\{1,r_b\}$.\\
2b. In case of death, choose $c^*\in\boldsymbol{c}^{(r,k-1)}$ uniformly, calculate the death ratio $r_d(\boldsymbol{c}^{(r,k-1)},c^{*})$ in \eqref{eq:death_ratio} and accept $\boldsymbol{c}^{(r,k)}=\boldsymbol{c}^{(r,k-1)}\setminus c^*$ with probability $\min\{1,r_b\}$\\
3. In case of non-acceptance in 2a. or 2b., set $\boldsymbol{c}^{(r,k)}=\boldsymbol{c}^{(r,k-1)}$\\
4. Calculate the probabilities $\pi^{\theta,k}_{i,j}, \pi^{\phi,k}_{i,j}$ as defined in \eqref{eq:piPhi} for $\boldsymbol{c}^{(r,k)}$ and $\boldsymbol{\beta}^{(r-1)}$\\
5. Set $k\rightarrow k+1\ls K$ and return to 1.

\textbf{II. M-step (Decoupled density estimation)}\\
6. Based on the MCMC-generated set of immigrants and probabilities, $\left(\boldsymbol{c}^{(r,k)},\pi^{\theta,k}_{i,j}, \pi^{\phi,k}_{i,j}\right)_{k=1}^K$, maximise the decoupled objective function \eqref{eq:obj_fct_detail} to get
\begin{equation*}
\boldsymbol{\beta}^{(r)}=\argmax_{\boldsymbol{\beta}} Q(\boldsymbol{\beta}\mid\boldsymbol{\beta}^{(r-1)})
\end{equation*}
7. Set $r\rightarrow r+1$ and return to step 1.
\end{mdframed} 
The algorithm may terminate after the parameters have converged, according to the selected criterion.

Modifications of the algorithm for the unmarked case, and for immigrants excluded, are relatively straightforward and discussed in Appendix \ref{sec:app3}.
  
\subsection{Comments}\label{sec:em_comments}
The algorithm requires storage and computation with matrices that are $O(n^2)$. On a standard PC, this makes computation prohibitive for samples with $n>10^4$. However this implementation is crude as, in this case, even the largest interevent time $t_n-t_1$ is considered as an interevent time by which $t_n$ could be triggered via $\theta(t_n-t_1)$ or $\phi(t_n-t_1)$, despite the fact that the probability of this could be effectively 0. Thus, for window size $t$ large relative to the support of $\theta$ and $\phi$, one can safely omit interevent times above a certain threshold. The result of this is banded/sparse matrices which reduce storage and computation from $n^2$ to $n \times m$, potentially with $m<<n$. In the case of the INARMA model, storage is less of an issue as only counts in bins need to be stored, rather than the location of each point.
  
  Technically, the use of MCMC in the E-step of the EM algorithm makes it a Monte-Carlo EM (MCEM) algorithm, having slightly weaker properties than the pure EM \cite{wei1990monte}. Regarding the convergence of the EM algorithm \cite{salakhutdinov2003optimization,Wu1983,Dempster1977}, one first needs that the necessary MLE regularity conditions are satisfied \cite{Ogata1978}, for instance having smooth distributions that are not too heavy tailed. Next, it must be ensured that the sequence of parameter estimates does not reach the boundary of the parameter space. For instance, if estimates $\mu$, $\eta$, or $\gamma$ are equal to zero at any iteration, or equivalently, if the support of $f(.)$ or $g(.)$ is smaller than the smallest interevent time, then the estimates will remain zero \eqref{eq:etaEMarma}. However, given non-zero starting estimates, the EM algorithm estimates satisfy the constraints of the model parameters. Regarding speed of convergence, there is the general result of \cite{Dempster1977,Wu1983} that the algorithm will not worsen the likelihood with each iteration. Further, from \cite{salakhutdinov2003optimization}, given that Q \eqref{eq:Q} is differentiable in $\boldsymbol{\beta}$ and the M-step has a unique solution, then the EM algorithm iterates in a positive direction on the true likelihood surface. Finally, when the missing information is small compared to the complete information, EM exhibits approximate Newton behavior with superlinear convergence near the true optimum. In terms of the ARMA point process, as well as other mixture type models, this means that, when clusters are overlapping, convergence will be slow, as has been shown for the Hawkes process with exponential offspring distribution \cite{LewisMohler2011,SornetteUtkin2009}, as well as other mixture models \cite{XuJordan1995,salakhutdinov2003optimization}. In particular for the NS process with immigrants included, since the number of offspring in a cluster is Poisson distributed, it will be difficult to distinguish heavily overlapping clusters from a single homogenous Poisson process. Exactly deriving the convergence properties of the ARMA point process model for a given parameterization would add no general insight.
     
  
\section{Simulation Study}\label{sec:simstudy}

Here a range of simulation studies are given to demonstrate performance and some relevant issues. Simple model specifications are employed: we use exponential kernels as in eq~\eqref{eq:renorm} with scale parameters $\theta_1^{-1}$ and $\phi_1^{-1}$, unmarked or exponentially distributed marks, moderate self-excitation, and constant immigration. Simulation is very fast, with details in appendix \ref{sec:app4}.

First is a demonstration of the performance of MCEM estimation of the simplest specification of the marked ARMA point process \eqref{eq:armamarks}:  exponential kernels with short memory, exponential mark distribution, moderate self-excitation, and constant immigration. The results of the repeated simulation and estimation are summarized in Tab~\ref{tab:sim1}. The main insight is that the estimation performs well. Further, as a robustness test, allowing for the immigration to be too flexible -- having up to 20 degrees of freedom, where the true immigration has only 1 -- has little impact on the estimated parameters. However, estimation is likely to be less robust to such an error when the immigration can better approximate the clusters (e.g., when clusters are larger and longer).

    \begin{table}[h!]
    \begin{center}
    \begin{tabular}{ c | c c c c c c }
        \toprule	
	p				& 0			& 3				& 5 				&	10 			&15 			&	20 \\
	\midrule
	$\mu$		& 0.979 	&	0.0999 	&	0.0995 	&	0.1010 	&	0.1014 	&	0.1019 \\
$\gamma$ 	& 4.21 	&	4.27 		&	4.27 		& 	4.16 		&	4.13  		&	4.03 \\
$\eta$   		&	0.51 	&	0.51 		&	0.50  		&	0.50  		&	0.51  		&	0.50 \\
    \bottomrule
    \end{tabular}
    \end{center}
        \caption{  Average parameter estimates of the marked ARMA, excluding memory scale parameters, averaged over 300 independent replications (simulations and estimations) done for a range of degrees of freedom (p) of the estimated immigration, using the R:logspline estimator. The standard deviation of the estimated parameters are about 0.008, 0.38, and 0.02 resp. (for all p). The simulated data is from the marked ARMA with parameters $\mu=0.1$, $\gamma=4$, $\theta_1^{-1}=0.1$, $\eta=0.5$, $\phi_1^{-1}=1$, with average sample size of 2000. In this case,  $\gamma$ is the scale parameter of the exponential mark distribution. Estimation is done by the MCEM algorithm, with K=50 MC samples taken from the chain of 300'000 Metropolis Hastings iterations.  }\label{tab:sim1}
    \end{table}

Next, acknowledging the lack of accessible likelihood as well as residuals to perform model selection or testing, the ability of the MCEM algorithm to ``select" the correct (sub-)model within the ARMA framework is important. E.g., on an AR simulation, the fitted ARMA should converge to have negligible $\gamma$ and consistently recover the AR part. To test this consistency, pure AR, pure MA, and the full ARMA are fit to AR, MA, and ARMA simulations, respectively, with results summarized in Tab.~\ref{tab:sim2}. For both AR and MA simulated data, the fitted ARMA converges well to the true sub-model. And for the full ARMA model, both AR and MA parts are well estimated. Of course, the range of misspecified models (e.g., the AR fit to the ARMA simulation) also provide estimates. Without an objective criterion to compare these different fits, it is therefore best to start with the broadest overall model -- in this case the ARMA -- and allow the MCEM algorithm to converge and select the relevant nested model. However, better methods to compare and test models are highly desirable. 

    \begin{table}[h!]
    \begin{center}
    \begin{tabular}{ c | c c c c c  }
        \toprule	
Fitted				& $\mu$		  & $\gamma$		& $\theta_1^{-1}$	& $\eta$			& $\phi_1^{-1}$		\\
	\midrule
ARMA sim.		&&&&& \\
	\midrule
ARMA 			& 0.095 (0.00025) 	& 5.52 (0.288) 	& 2.03 (0.12) 		& 0.47 (0.0024) 	& 0.097 (0.00019) 	\\
NS/MA			& 0.21  (0.004)   		& 4.84 (1.62)  	& 0.67 (0.049) 	& --							& --							\\
Hawkes/AR	& 0.27  (0.003)   		& --						& -- 						& 0.77 (0.0014)		& 0.22 (0.0015)		\\
	\midrule
MA sim.			&&&&& \\
	\midrule
ARMA 			& 0.10 (0.00016)  	& 6.98 (0.10) 		& 1.96 (0.017) 	& 0.01 (0.0001) 	& 2.06 (5.55) 	\\
NS/MA			& 0.10  (0.002)   	& 7.06 (0.09)  	& 1.98 (0.016) 	& --							& --						\\
Hawkes/AR	& 0.16  (0.001)  	& --						& -- 						& 0.81 (0.0007)		& 0.83 (0.0066)		\\
	\midrule
AR sim.			&&&&& \\
	\midrule
ARMA 			& 0.09 (0.0002)   		& 0.13 (0.10) 		& 5.58 (54.2) 		& 0.79 (0.003) 		& 0.099 (0.0001) 	\\
NS/MA			& 0.13  (0.002)   		& 2.73 (0.33)  	& 0.18 (0.0006)	& --							& --							\\
Hawkes/AR	& 0.10  (0.0001)  		& --						& -- 						& 0.79 (0.0023)		& 0.10 (0.00001)	\\
    \bottomrule
    \end{tabular}
    \end{center}
        \caption{ Average (and standard deviation) of  parameter estimates over 100 independent replications. The first set are for data simulated from the ARMA with parameters $(\mu=0.1,\gamma=5,\theta_1^{-1}=2,\eta=0.5,\phi_1^{-1}=0.1)$ with the same unmarked ARMA \eqref{eq:armaP}, and nested MA, and AR sub-models fit respectively. The second set are for data simulated from the MA with parameters $(\mu=0.1,\gamma=7,\theta_1^{-1}=2)$, and the final set for data simulated from the AR (Hawkes) with parameters $(\mu=0.1,\eta=0.8,\phi_1^{-1}=0.1)$.  All simulations contain 500 points, and the MCEM estimation is done with K=50 MC samples taken from a chain of length 300'000 generated by the Metropolis-Hastings algorithm.  
 }\label{tab:sim2}
    \end{table}

Finally, the effect of cluster overlap is briefly examined, by fitting the unmarked ARMA, with exponential kernels, and varying the immigration intensity. It is intuitive that more heavily overlapping clusters are more difficult to distinguish from immigration. As summarized in Table~\ref{tab:sim3}, as the clusters becomes more overlapping, indeed the immigration intensity becomes increasingly overestimated and in particular $\gamma$ underestimated.

    \begin{table}[h!]
    \begin{center}
    \begin{tabular}{ c | c c c c c c c c c c  }
    \toprule	
$\mu_{sim}$ 				& 0.10 & 0.31 & 0.52 & 0.73 & 0.94 & 1.16 & 1.37 & 1.58 & 1.79 & 2.00 \\
	\midrule
$\mu$ 					  		& 0.10 & 0.32 & 0.57 & 0.78 & 1.05 & 1.33 & 1.62 & 1.80 & 2.04 & 2.23 \\
$\gamma$					& 6.17 & 5.97 & 5.76 & 5.19 & 4.76 & 4.36 & 4.21 & 4.23 & 3.89 & 3.74 \\
$\theta_1^{-1}$ 		& 0.54 & 0.54 & 0.52 & 0.47 & 0.44 & 0.42 & 0.40 & 0.40 & 0.37 & 0.34 \\
$\eta$ 						& 0.39 & 0.39 & 0.40 & 0.47 & 0.47 & 0.49 & 0.49 & 0.54 & 0.55 & 0.60 \\
$\phi_1^{-1}$					& 0.24 & 0.24 & 0.25 & 0.28 & 0.29 & 0.31 & 0.31 & 0.31 & 0.31 & 0.30 \\
sd($\mu$) 					& 0.011 & 0.041 & 0.085 & 0.165 & 0.258 & 0.342 & 0.435 & 0.542 & 0.605 & 0.773 \\
sd($\gamma$) 			& 0.426 & 0.740 & 0.956 & 1.288 & 0.994 & 1.284 & 1.091 & 1.350 & 1.103 & 1.734 \\
sd($\theta_1^{-1}$) 	& 0.041 & 0.060 & 0.072 & 0.095 & 0.092 & 0.097 & 0.075 & 0.078 & 0.052 & 0.066 \\
sd($\eta$) 					& 0.059 & 0.096 & 0.107 & 0.115 & 0.105 & 0.132 & 0.113 & 0.126 & 0.137 & 0.134 \\
sd($\phi_1^{-1}$)	 		& 0.033 & 0.040 & 0.042 & 0.071 & 0.078 & 0.104 & 0.072 & 0.068 & 0.072 & 0.064 \\
    \bottomrule
    \end{tabular}
    \end{center}
        \caption{  Average (top rows) and standard deviation (lower rows) for estimated parameters, computed over 500 independent replications for each value of immigration intensity used in the simulation, $\mu_{sim}$. The data are simulated from the unmarked ARMA with parameters $(\mu_{sim},\gamma=5,\theta_1^{-1}=0.5, \eta=0.5, \phi_1^{-1}=0.25)$ and with each realization containing 1000 points. The MCEM estimation is done with K=50 MC samples taken from a chain of length 300'000 generated by the Metropolis-Hastings algorithm.  
 }\label{tab:sim3}
    \end{table}

   \pagebreak
   \clearpage
   
\section{Appendix} 
\subsection{Derivation of the palm intensity in the exponential case}\label{sec:app1}
Below we will derive an explicit solution for the palm intensity function \eqref{eq:PI} of the ARMA point process for exponential densities. In this goal, let the setting in section \ref{sec:def} be fulfilled, with $\theta_0,\theta_1,\phi_0\in(0,\infty)$, $\phi_1\in(\phi_0,\infty)$, $\theta(t)=\theta_0 e^{-\theta_1 t}$, $\phi(t)=\phi_0 e^{-\phi_1 t}$, with $h$ and $\tau$ be defined by $\eqref{eq:PI}$ and $\eqref{eq:tau1}$, respectively, and $\rho=\overline{\lambda}h$. The Laplace transform of a function $f:[0,\infty)\rightarrow\R$ will be denoted by $f^\ast$.

To get an explicit solution for $\tau$, we Laplace-transform equation \eqref{eq:tau2} to get
\begin{equation}
\tau^{\ast}(s)=\dfrac{\mu^2(1+\gamma)}{s}+\dfrac{\mu\phi_0+\overline{\lambda}\mu\eta}{\phi_1+s}+\dfrac{\mu\theta_0}{\theta_1+s}+\dfrac{\phi_0}{\phi_1+s}\tau^{\ast}(s).
\end{equation}
Rearranging terms to solve for $\tau^{\ast}(s)$ and doing a partial fraction decomposition allows us to invert the Laplace transform to arrive at 
\begin{as}
\tau(u)&=\overline{\lambda}\mu+\left(L_1 e^{-\left(\phi_1-\phi_0\right)u}+L_2 e^{-\theta_1 u}\right)\bm{1}\{u>0\}, \text{ where }\\
L_1&=\mu\phi_0\dfrac{\theta_1+\theta_0-\phi_1+\phi_0}{\theta_1-\phi_1+\phi_0},\\ L_2&=\mu\theta_0\dfrac{\theta_1-\phi_1}{\theta_1-\phi_1+\phi_0}.
\end{as}
Using this explicit solution for $\tau$, we get from equation \eqref{eq:PIsol} that
\begin{as}\label{eq:aux}
\rho(u)&=\mu\overline{\lambda}(1+\gamma)+e^{-\theta_1 u}\left[\mu\theta_0+\dfrac{L_1\theta_0}{\theta_1+\phi_1-\phi_0}+\dfrac{L_2\gamma}{2}\right]+e^{-\phi_1 u}\overline{\lambda}\phi_0\\
&\quad+\int_0^u\phi(x)\rho(u-x)dx+\int_0^\infty\phi(x+u)\rho(x)dx.
\end{as}
Analogously to the solution for $\tau$, we do a Laplace transformation of equation \eqref{eq:aux}, partial fraction fraction decomposition and an inverse transform to arrive at
\begin{as}
\rho(u)&=\overline{\lambda}^2+K_1 e^{-\left(\phi_1-\phi_0\right)u}+K_2 e^{-\theta_1 u},
\end{as}
where the constants $K_1,K_2\in[0,\infty)$ are given by
\begin{as}
K_0&=\mu\theta_0+\mu\phi_0\dfrac{\theta_1+\theta_0-\phi_1+\phi_0}{\theta_1-\phi_1+\phi_0}\dfrac{\theta_0}{\theta_1+\phi_1-\phi_0}+\dfrac{1}{2}\dfrac{\mu\theta_0(\theta_1-\phi_1)\theta_0}{(\theta_1-\phi_1+\phi_0)\theta_1}\\
K_1&=K_0
\left(\dfrac{\phi_1\phi_0}{\theta_1+\phi_1}-\dfrac{\phi_0}
{\phi_1-\phi_0-\theta_1}\right)+\overline{\lambda}\phi_0\left(1+\dfrac{\phi_0}{2(\phi_1-\phi_0)}\right)\\
K_2&=K_0\dfrac{\phi_1-\theta_1}{\phi_1-\phi_0-\theta_1}.
\end{as}

\subsection{EM algorithms for the ARMA point process}\label{sec:app3}

\subsubsection{Unmarked ARMA}
If $T$ is an ARMA point process as in section \ref{sec:deriv} without marks and $\sigma(N^\mu,N)$-conditional intensity
\begin{equation}
\lambda(s)
= \mu(s) + \sum_{j=-\infty}^{N^{\mu}(s)} \theta(s-T^{\mu}_j) + \sum_{k=-\infty}^{N(s)} \phi(s-T_k)~,
\end{equation}
estimation is almost identical. The main difference is that $\gamma=\int_0^\infty \theta(s) ds$ will be included as a parameter instead of the density of marks, $m$. The complete data likelihood \eqref{eq:likIntPart} becomes 
\begin{as}
    &L(\boldsymbol{\beta}\mid \boldsymbol{t},\boldsymbol{c},\boldsymbol{Z}) = \prod_{i=1}^{n_c} \mu(s_i) \text{Exp}\lbrace - \int_0^t \mu(s)ds \rbrace \times  \\
    & ~~~\prod_{i=1}^{n}\prod_{j=1}^{n_c} \Big[ \theta(t_i-s_j) \Big]  ^{ Z_{i,j}^\theta } \text{Exp}\lbrace -   \sum_{j=1}^{n_c} \int_0^t \theta(s-s_j) ds \rbrace \times \\
    & ~~~~~~\prod_{i=1}^{n}\prod_{k=1}^{n} \Big[ \phi(t_i-t_k) \Big] ^{Z_{i,k}^\phi} \text{Exp}\lbrace - \sum_{j=1}^{n}\int_0^t \phi(s-t_j)ds \rbrace~
\end{as}
and $\gamma$ is estimated by
\begin{equation}
\widehat{\gamma}=\frac{ \sum_{i,j} \pi^\theta_{i,j} }{ \sum_{j=1}^{n_c} \int_0^{t-s_{j}}g(s)ds }.
\end{equation}
All equations of sections \ref{sec:immknown}, \ref{sec:immunknown} and \ref{sec:em_stepwise} extend using constant marks $y\equiv 1$.

\subsubsection{Immigrants not included in the sample}
As introduced in section \ref{sec:comments}, one may wish to estimate a modification of the ARMA point process, where the immigrant process $N^{\mu}$ is not included in the process $N$ and the $\sigma(N^\mu,N)$-conditional intensity is given by
\begin{equation}
\lambda(s)
= \mu(s) + \sum_{j=-\infty}^{N^{\mu}(s)} Y_j \theta(s-T^{\mu}_j) + \sum_{k=-\infty}^{N(s)} \phi(s-T_k)~.
\end{equation}
To describe the (mostly straightforward) extension of the algorithm above, additional to the notation of section \ref{sec:deriv}, we will need $\boldsymbol{\tilde{t}}=\boldsymbol{t}\cup\boldsymbol{s}$, the \emph{ordered} union of points and immigrants with $\tilde{n}=n+n_c$ elements. The indicator variables \eqref{eq:Z} have a modified domain
\begin{equation}
  \boldsymbol{Z}=\lbrace {Z}^{\theta}_{i,j}~,~i=1,...,n,~j=1,...,n_c \rbrace \cup \lbrace {Z}^{\phi}_{i,j}~,~i=1,...,\tilde{n},~j=1,...,\tilde{n} \rbrace~,
\end{equation}
and the complete data likelihood becomes
\begin{as}
&L(\boldsymbol{\beta}\mid \boldsymbol{t},\boldsymbol{c},\boldsymbol{Z}) =     
\prod_{i=1}^{n_c} \mu(s_i) \text{Exp}\lbrace - \int_0^t \mu(s)ds \rbrace m(y_i) \times \nonumber \\
    & ~~~\prod_{i=1}^{n}\prod_{j=1}^{n_c} \Big[ y_j \theta(t_i-s_j) \Big]  ^{ Z_{i,j}^\theta } \text{Exp}\lbrace -   \sum_{j=1}^{n_c} y_j \int_0^t \theta(s-s_j) ds \rbrace \times \nonumber \\
    & ~~~~~~\prod_{i=1}^{n}\prod_{k=1}^{\tilde{n}} \Big[ \phi(t_i-\tilde{t}_k) \Big] ^{Z_{i,k}^\phi} \text{Exp}\lbrace - \sum_{j=1}^{\tilde{n}}\int_0^t \phi(s-\tilde{t}_j)ds \rbrace~.
\end{as}

Moreover, the MCMC algorithm used in the E-step needs to be adapted. Given the current state of the Markov chain $\boldsymbol{c}=(\boldsymbol{s},\boldsymbol{y})$, in a birth step, we choose an immigrant $s^*$ uniformly in the time window $[0,t]$ (instead of choosing among offspring), while a mark $y^*$ is generated from $m$. As before, in a death step, we choose an immigrant uniformly among the immigrant set $\boldsymbol{c}$. The Metropolis-Hastings birth ratio becomes 
\begin{as}
r_b(\boldsymbol{c},c^{*})=
&\frac{\vert t\vert}{n_c+1} \mu(s^{*})\text{Exp}\lbrace -\int_{0}^{t} y^{*}\theta(s-s^{*}) +\phi(s-s^{*})ds \rbrace \times
\\ 
&\prod_{t_i\in \boldsymbol{t}} \left(  1+ \frac{y^{*}\theta(t_i-s^{*})+\phi(t_i-s^{*}) }{ \sum_{s_j\in \boldsymbol{s}} y_j\theta(t_i-s_j)+\sum_{ \tilde{t}_j\in \boldsymbol{\tilde{t}} } \phi(t_i-\tilde{t}_j) }  \right)~
\end{as}
and the death ratio $r_d(\boldsymbol{c},c^{*})=1/{r_b(\boldsymbol{c}\setminus c^{*},c^*)}$.

Then the algorithm described in section \ref{sec:em_stepwise} carries over with the obvious modifications of probabilities \eqref{eq:piPhi} and objective function \eqref{eq:obj_fct_detail}.

  \subsection{Simulation of the ARMA point process}\label{sec:app4}
   
Below a detailed algorithm for simulating ARMA point processes is presented. It exploits the fact that, given the branching structure, the innovation, MA, and AR processes are mutually independent inhomogeneous Poisson processes. Efficient simulation algorithms for the Hawkes process follow the same approach.  At the end of step II, one has simulated an NS process. By skipping step II, and completing step III, one simulates a Hawkes process. To avoid edge effects, one should simulate on a large window, and discard the \emph{burn in period}. When the kernels can be integrated and inverted, using inverse transform sampling \cite{cinlar2013introduction} makes the algorithm very fast. 

 \subsubsection*{Simulation algorithm}

 \setlength{\parindent}{0cm}
\begin{mdframed}[leftmargin=0pt,rightmargin=0pt]

  \textbf{I. Simulate the immigrant points}\\
  (i) Simulate Poisson process $\lbrace T_i^{(0)} \rbrace_{i\in1,\dots,n_c}$ on the window $(0,t]$, where $N(t)=n_c$.\\
  \textbf{II. Simulate MA points}\\
  (i) For each immigrant $i=1,\dots,n_c$: simulate the number of offspring $N_i^\theta \overset{i.i.d}{\sim}\text{Pois}(\gamma)$, and then sample the $N_i^\theta$ inter-event times $S_{i,j},~j=1,\dots,N_i^\theta$, i.i.d from pdf $g$. If $N_i^\theta=0$, simulate zero inter-event times for that immigrant.\\
  (ii) The MA points generated by the $i^{th}$ immigrant are then $\lbrace T^{\theta}_{i,j} \rbrace_{j=0:N_i^\theta}=T^{(0)}_i+\lbrace S_{i,j} \rbrace_{j=0:N_i^\theta}$. \\
  (iii) The immigrant and MA points together are $\lbrace T_i^{(1)} \rbrace_{i\in1,\dots,n_1}=\lbrace T_i^{(0)} \rbrace \cup \lbrace T^{\theta}_{1,j} \rbrace \cup \dots \cup \lbrace T^{\theta}_{n_c,j} \rbrace$, where $n_1=n_0+\sum_{i=1}^{n_c} N_i^\theta$.\\
  \textbf{III. Simulate AR points by generation}\\
  (i) Set the fertile points $A=\lbrace 1,\dots,n_1 \rbrace$, generation $k=1$, and zeroeth generation points $\lbrace T^{\phi[0]}_i \rbrace=\lbrace T_i^{(1)} \rbrace$. \\ 
  (ii) For the current generation k, for all fertile points $\forall i\in A$: simulate the number of direct offspring $N_i^{\phi[k]}\overset{i.i.d}{\sim}\text{Pois}(\eta)$, and then the $N_i^{\phi[k]}$ inter-event times $S^{\phi[k]}_{i,j},~j=1,\dots,N_i^{\phi[k]}$, i.i.d from pdf $f$.  If $N_i^{\phi[k]}=0$, simulate zero inter-event times for that point.\\
  (iii) The AR points generated by the $i^{th}$ point are then $\lbrace T^{\phi[k]}_{i,j} \rbrace_{j=0:N_i^{\phi[k]}}=T^{\phi[k-1]}_i+\lbrace S^{\phi[k]}_{i,j} \rbrace_{j=0:N_i^{\phi[k]}}$ and \\
  (iv) the union of these sets, $\lbrace T^{\phi[k]}\rbrace = \lbrace T^{\phi[k]}_{1,j} \rbrace \cup \dots \cup \lbrace T^{\phi[k]}_{N_i^{\phi[k]},j} \rbrace$, is the offspring of generation k. 
  (v) Update the fertile set $A=\lbrace i: T_i^{\phi[k]} < r \rbrace$ to be all points born in the current generation k that fall within the window $(0,t]$.\\
  (vi) If $A$ is non-empty, then increment the generation ($k=k+1$) and return to (ii), otherwise return the realization formed by joining all generations:
   $\{T_i\}_{i=1:n}=\lbrace T_i^{(1)} \rbrace \cup \lbrace T^{\phi[1]}\rbrace \cup \dots \cup \lbrace T^{\phi[k]}\rbrace $.\\

\end{mdframed}
 

\cleardoublepage
\phantomsection
\bibliographystyle{abbrv}
{\footnotesize
\bibliography{ARMA}
}

\end{document}